\documentclass[11pt]{article}
\usepackage{amssymb}\usepackage{amsmath}
\newcommand{\fr}[2]{{\textstyle \frac{#1}{#2} }}
\newcommand{\be}{\begin{eqnarray}}
\newcommand{\ee}{\end{eqnarray}}
\newcommand{\nn}{\nonumber}
\newcommand{\tr}{\mathop{\rm tr}\nolimits}
\newcommand{\one}{{\mathfrak 1}}
\newcommand{\two}{{\mathfrak 2}}
\newcommand{\three}{{\mathfrak 3}}
\newcommand{\et}{\eta_{\scriptscriptstyle 0}}
\newcommand{\proof}{{\em Proof.\ }}

\def\qed{\hfill $\square$}
\newcommand{\+}{{+}}
\renewcommand{\=}{{=}} 
\renewcommand{\-}{{-}} 
\newcommand{\FRW}[6]{\Bigl\{\begin{smallmatrix}%
 \textstyle #1\, & \textstyle #3\, & \textstyle #5 \\[2pt]%
 \textstyle #2\, & \textstyle #4\, & \textstyle #6 %
 \end{smallmatrix}\Bigr\}}
\newcommand{\FCG}[6]{\Bigl[\begin{smallmatrix}%
 \textstyle #1\, & \textstyle #3\, & \textstyle #5 \\[2pt]%
 \textstyle #2\, & \textstyle #4\, & \textstyle #6 %
 \end{smallmatrix}\Bigr]}
\newtheorem{prop}{Proposition}
\newtheorem{lem}{Lemma}
\newcounter{num}
\newcommand{\rem}{{\em Remark~\thenum.}\ \addtocounter{num}{1}}
\begin{document}
December 2004  \strut\hfill math.QA/0412482 \\
\vspace*{3mm}
\begin{center}
{\large\bf
On higher spin $U_q(sl_2)$--invariant $R$--matrices
} \\ [3mm]
{\sc Andrei G. Bytsko} \\ [2mm]
{ Steklov Mathematics Institute \\ 
 Fontanka 27, 191023 St.~Petersburg, Russia } \\ [3mm]
{Dedicated to Professor L.D.~Faddeev on the occasion 
  of his 70th birthday }

\end{center}
\vspace{1mm}
\begin{abstract}

The spectral decomposition of regular $U_q(sl_2)$--invariant
solutions of the Yang--Baxter equation is studied.
An algorithm for finding all possible solutions of 
spin~$s$ is developed. It also allows to reconstruct 
the R--matrix {}from a given nearest neighbour spin chain 
Hamiltonian. The algorithm is based on reduction of
the Yang--Baxter equation to certain subspaces. 
As an application, the complete list of inequivalent
regular $U_q(sl_2)$--invariant R--matrices is obtained
for generic~$q$ and spins $s\leq\frac{3}{2}$.
Some further results about spectral decompositions for
higher spins are also obtained. In particular, it is
proved that certain types of regular $sl_2$--invariant 
R--matrices have no $U_q(sl_2)$--invariant counterparts.

\end{abstract}
\section{Preliminaries}   

The quantum Lie algebra $U_q(sl_2)$ is defined as a universal 
enveloping algebra over $\mathbb C$ with the identity element
$e$ and generators $S^\pm$, $S^z$ obeying the following 
defining relations~\cite{KuR}
\begin{equation}\label{cS}
  [S^+ , S^- ] = [ 2 S^z ]_q \,,\qquad
  [S^z , S^\pm  ] = \pm S^\pm \,,
\end{equation}
where
$ [t]_q = (q^t \- q^{-t})/(q \- q^{-1}) $. 
$U_q(sl_2)$ can be equipped with a Hopf algebra
structure~\cite{SkDr}. In particular, the 
co--multiplication (a co--associative linear homomorphism)
is defined as follows
\begin{equation}\label{del}
 \Delta(S^\pm) = S^\pm \otimes q^{-S^z} + q^{S^z} \otimes S^\pm \,,
 \qquad \Delta(S^z) = S^z \otimes e + e \otimes S^z \,.
\end{equation}
For generic~$q$, the algebra (\ref{cS})--(\ref{del}) has the same 
structure of representations as the undeformed algebra 
$sl_2$~\cite{Ro}. In particular, irreducible 
highest weight representations $V_s$ are parameterized by 
a non--negative integer or half--integer number $s$
(referred to as {\em spin}) and are $(2s{+}1)$--dimensional.
We will use the standard notation, $|k\rangle$, $k=-s\,\ldots,s$,  
for the basis vectors of~$V_s$ such that 
$S^z |k\rangle = k |k\rangle$, $\langle k'|k\rangle=\delta_{kk'}$.

Let $\mathbb E$ denote the identity operator 
on~$V_{s}^{\otimes 2}$.
Consider operator valued functions,
$R(\lambda) : {\mathbb C} \mapsto {\rm End}\ V_{s}^{\otimes 2}$,
that have the following properties:
\be
\label{reg}
 regularity: && R(0) = \mathbb E \,, \\ [1mm]
\label{uni}
 unitarity: &&
 R(\lambda)\, R(-\lambda) = \mathbb E \,,\\ 
\label{ans1}
 spectral\ decomposition:
 && R(\lambda) = \sum_{j=0}^{2s} 
 r_j(\lambda) P^{j} \,, \\ 
 normalization: && \label{gauge}
 r_{2s}(\lambda) = 1 \,.
\ee
Here $P^j$ is the projector onto the spin $j$ 
subspace $V_j$ in~$V_s^{\otimes 2}$
and $r_j(\lambda)$ is a scalar function.
Property (\ref{ans1}) is equivalent to the 
requirement of $U_q(sl_2)$--{\em invariance}, i.e.,
\be
 [ R(\lambda) , \Delta(\xi) ] = 0 \qquad
 \forall \, \xi \in U_q(sl_2) \,.
\ee
In order to fulfill (\ref{reg})--(\ref{uni}), the 
coefficients $r_j(\lambda)$ must satisfy the relations
\be
\label{ans2}
 && r_j(0) = 1 \,, \qquad r_j(\lambda) r_j(-\lambda) =1 \,.
\ee
In what follows we will assume that $r_j(\lambda)$ 
are analytic in some neighbourhood of~$\lambda=0$. 
The normalization condition (\ref{gauge}) is imposed 
in order to eliminate inessential freedom of rescaling 
$R(\lambda)$ by an arbitrary analytic function 
preserving the conditions (\ref{ans2}).

For a given $R(\lambda)$ satisfying (\ref{reg})--(\ref{uni}), 
let us define the Yang--Baxter (YB) {\em operator},
$Y(\lambda,\mu): {\mathbb C}^2 \mapsto 
 {\rm End}\ V_{s}^{\otimes 3}$, as follows
\be\label{YBO}
 Y(\lambda,\mu) =
 R_{\one\two}(\lambda)\, R_{\two\three}(\lambda+\mu) \, 
 R_{\one\two}(\mu) -  R_{\two\three}(\mu)\, 
 R_{\one\two}(\lambda+\mu)\, R_{\two\three}(\lambda) \,. 
\ee
Here and below we use the standard notations -- the lower 
indices specify the components of the tensor 
product~$V_{s}^{\otimes 3}$. We will say that $R(\lambda)$ is 
an ($U_q(sl_2)$--invariant) {\em R--matrix} 
if the corresponding YB operator vanishes
on~$V_{s}^{\otimes 3}$,
\be\label{YBE}
 Y(\lambda,\mu) = 0 \,. 
\ee
An advantage of treating the YB {\em equation} (\ref{YBE}) 
as the condition of vanishing of the YB operator is that, 
as we will show below, conditions of vanishing of the YB 
operator on some {\em subspaces} of~$V_{s}^{\otimes 3}$ 
involve fewer coefficients~$r_j(\lambda)$. 
Moreover, $r_j(\lambda)$ found by resolving such
condition for a given subspace can be further used 
in order to write down and solve conditions of vanishing 
of the YB operator on other subspaces. A recursive procedure
of this type will be presented in the next section.

\rem Eqs.~(\ref{ans2}) and (\ref{YBE}) are preserved under 
rescaling of the spectral parameter,
\be\label{ren}
  \lambda \rightarrow \gamma \, \lambda \,,
\ee
by an arbitrary finite constant~$\gamma$. We will regard 
R--matrices related by such transformation with a finite 
nonzero $\gamma$ as equivalent.

\rem 
Conditions (\ref{reg}) and (\ref{gauge}) along with the YB 
equation ensure unitarity of an R--matrix (see Appendix~A).

For $q=1$, there are known~\cite{RM,KRS,ZZ,Ba}
four different types of $sl_2$--invariant R--matrices:
\be
\label{ser0a}
 && R(\lambda) = (1- \lambda)^{-1} \Bigl(
  \mathbb E - \lambda \  \mathbb P \Bigr) \,, \\ [1mm]
\label{ser0b}
 && R(\lambda) = P^{2s} +
 \sum_{j=0}^{2s-1} \Bigl( \prod_{k=j+1}^{2s}
 \frac{k+\lambda}{k-\lambda}  \Bigr) \,  P^{j} \,,\\
\label{ser0c}
 && R(\lambda) = (1 - \lambda)^{-1}
 \Bigl( \mathbb E -  \lambda \mathbb P + 
     \frac{\beta\lambda}{\lambda - \alpha} P^0 \Bigr) \,,  \\
\nn 
 && \qquad
      \qquad \alpha= s + \fr{1}{2} + (-1)^{2s+1} \,, \quad
      \beta = (2s+1)(-1)^{2s+1} \,, \\
\label{ser0d}
 && R(\lambda) = 
 \mathbb E + (b^2+1) \frac{1-e^\lambda}{e^\lambda-b^2} P^0 \,,
 \qquad b+b^{-1}=2s+1 \,,
\ee
where 
\be\label{perm}
 \mathbb P = \sum_{j=0}^{2s} (-1)^{2s-j} P^j 
\ee
is the permutation in $V_{s} \otimes V_{s}$. 
Observe that for all but the last type of solutions we have 
\be\label{as0}
 R(\pm\infty)  \, = \, \mathbb P \,.
\ee 

For $s=\frac{1}{2}$, R--matrices (\ref{ser0b}) and (\ref{ser0c})
degenerate into (\ref{ser0a}) and the fourth solution, 
(\ref{ser0d}), is absent. For $s=1$,  R--matrices (\ref{ser0b}) 
and (\ref{ser0c}) are equivalent. For $q=1$ and $s=3$, there is 
known an additional solution which is not of the form 
(\ref{ser0a})--(\ref{ser0d}). It is given by~\cite{Ke}
\be\label{s3a}
 R(\lambda) = P^6 +  
 \frac{1 \+ \lambda}{1 \- \lambda } \, P^5 + P^4 +
  \frac{4 \+ \lambda}{4 \- \lambda } \, P^3 + P^2 +
 \frac{1 \+ \lambda}{1 \- \lambda } \, P^1 + 
 \frac{1 \+ \lambda}{1 \- \lambda } \, 
 \frac{6 \+ \lambda}{6 \- \lambda } \, P^0 \,.
 \ee
Numerical, computer--based checks~\cite{Ke} suggest that 
eqs.~(\ref{ser0a})--(\ref{ser0d}) and (\ref{s3a}) 
exhaust the list of $sl_2$--invariant R--matrices,
but no corresponding theorem has been proved yet.

The $q\neq 1$ counterparts of (\ref{ser0b}) and
(\ref{ser0d}) are given by \cite{Jim,Ba}
\be\label{serB}
 && R(\lambda) =  P^{2s} +
 \sum_{j=0}^{2s-1} \Bigl( \prod_{k=j+1}^{2s}
 \frac{[k+\lambda]_q}{[k-\lambda]_q}  \Bigr) \,  P^{j} \,, \\[0.5mm]
\label{qTL}
 && R(\lambda) = 
 \mathbb E + (b^2+1) \frac{1-e^\lambda}{e^\lambda-b^2} P^0 \,,
 \qquad b+b^{-1}=[2s+1]_q \,.
\ee

The aim of the present paper is to study 
$U_q(sl_2)$--invariant solutions of the Yang--Baxter 
equation for a generic $q$ (i.e., 
$q$ is not a root of unity and $q\neq 0,\infty$)
and to develop a systematic method of finding
all possible sets of $r_j(\lambda)$ for a given spin~$s$.
In particular, we will prove that (\ref{ser0a}),
(\ref{ser0c}), and (\ref{s3a}) have no regular 
$U_q(sl_2)$--invariant counterparts.
Our approach will be based on the fact that 
$U_q(sl_2)$--invariance of an \hbox{R--matrix} implies 
that the corresponding YB operator commutes with
the action of $U_q(sl_2)$ on~$V_s^{\otimes 3}$.
This action is defined as  
\be
\label{genSz}
 S^z_{\one\two\three} &=& (\Delta \otimes id) \, \Delta S^z
 = S^z_\one + S^z_\two + S^z_\three \,,\\ [1mm]
\label{genSpm}
 S^\pm_{\one\two\three} &=& (\Delta \otimes id) \, \Delta S^\pm
 = S^\pm_{\one\two} \, q^{-S^z_\three} + q^{S^z_{\one\two}} 
  \, S^\pm_{\three}  = q^{S^z_{\one}} \, S^\pm_{\two\three} 
   + S^\pm_{\one} \, q^{-S^z_{\two\three}} \,.
\ee
The assertion 
\be\label{YSS}
 [Y(\lambda,\mu), S^z_{\one\two\three}]= 
  [Y(\lambda,\mu), S^\pm_{\one\two\three}] =0
\ee
follows {}from the fact that $P^j$ are functions of 
$\Delta C$, where $C$ is the Casimir element of $U_q(sl_2)$.
It is obvious {}from (\ref{genSz})--(\ref{genSpm}) that
$P^j_l$, \hbox{$l={\scriptstyle \{12\}},{\scriptstyle \{23\} }$}
commute with $S^z_{\one\two\three}$ and~$S^\pm_{\one\two\three}$.

\section{Reduced Yang--Baxter equations}
\subsection{Hecke--Temperley--Lieb algebra in YB}

Interrelations between Hecke algebras, braid groups, and 
constant (independent on the spectral parameter~$\lambda$) 
solutions of the YB equation are well known. In the case 
of nontrivial spectral parameter dependence, a construction
of an R--matrix employing the Temperley--Lieb algebra \cite{TL}
generators was introduced by Baxter in~\cite{Ba}.
For the purposes of the present article, we will need 
the following, slightly generalized, version of this
construction. 

\begin{lem}\label{Hecke}
Consider an associative algebra over $\mathbb C$
with the unit element ${\mathbb E}$
and the generators $U_l$ labeled by
\hbox{$l={\scriptstyle \{12\}},{\scriptstyle \{23\} }$}
obeying the following Hecke--type
relations (with $\et$ and $\eta_{\scriptscriptstyle 1}$ 
being scalar constants, $\Re \et \geq 0$)
\be\label{Hecke1}
 &&  U_l^2 = \et \, U_l + \eta_{\scriptscriptstyle 1} 
 \, {\mathbb E} \,, \\
\label{Hecke2}
 && U_{\one\two} \, U_{\two\three} \, U_{\one\two} - 
 U_{\two\three} \, U_{\one\two} \, U_{\two\three} =
 U_{\one\two} - U_{\two\three} \,.
\ee
Let $g(\lambda)$ be a function analytic in
a neighbourhood of~$\lambda=0$ and satisfying 
the condition $g(0)=0$. 
Then \hbox{$R_l(\lambda)=\mathbb E + g(\lambda) \, U_l$}
satisfy the YB equation (\ref{YBE}) if and only if
\be\label{gsol}
 g(\lambda)= \left\{ \begin{array}{lll}
 \frac{2\gamma\lambda}{1-\gamma\lambda}  & &
  {\rm if}\ \et = 2 \,; \\  [1mm]
 \textstyle
  b\, \frac{1-e^{\gamma\lambda}}%
  {e^{\gamma\lambda} - b^2} \,, &
 \quad b+b^{-1}=\et &  {\rm if}\ \et \neq 2 \,.
 \end{array} \right. 
\ee
Here $\gamma$ is an arbitrary finite constant.
\end{lem}

\rem R--matrices
(\ref{ser0a}) and (\ref{ser0d}) provide two examples
where $U_l$ are elements of ${\rm End}\, V_s^{\otimes 3}$ 
given by 
\be\label{Ue}
 U_{\one\two}=U\otimes e, \qquad U_{\two\three}=e\otimes U 
\ee
with $U=\mathbb{ E + P}$, $\et=2$ and $U=P^0$, $\et=2s{+}1$, 
respectively, and $e$ being the unit element on~$V_s$.
Let us however stress that, in general, the hypotheses 
of Lemma~~\ref{Hecke} {\em do not} require that 
$U_l$ be of the form (\ref{Ue}) 
with $U \in {\rm End}\, V^{\otimes 2}$. 
In fact, below we will apply Lemma~\ref{Hecke} in the 
cases where the underlying linear space is not of 
the form~$V^{\otimes 3}$.

\proof For the sake of completeness of the exposition,
let us give a proof of the lemma. 
Substituting $R_l(\lambda)$ into (\ref{YBE}) and
employing (\ref{Hecke1})--(\ref{Hecke2}), one reduces
the YB equation to the form
$(...)(U_{\one\two} - U_{\two\three})=0$, where $(...)$ is a
scalar factor. Therefore, the YB equation holds if and only 
if this factor vanishes, which is equivalent to the requirement 
that $g(\lambda)$ satisfy the following functional relation
\be\label{gHecke}
 g(\lambda-\mu) \, g(\lambda) \, g(\mu) + 
 \et \, g(\lambda-\mu) \, g(\mu) + 
 g(\lambda-\mu) - g(\lambda) + g(\mu) = 0 \,.
\ee 
By differentiating (\ref{gHecke}) w.r.t.~$\mu$, setting 
$\mu=\lambda$, and taking into account the condition $g(0)=0$, 
one derives the differential equation
\be\label{gdif}
 g^\prime(\lambda) =
 g^\prime(0) \, \bigl( (g(\lambda))^2 + \et \, g(\lambda) 
  + 1 \bigr) \,.
\ee
Its solution is given by (\ref{gsol}) 
and it is easily verified that this solution does 
satisfy~(\ref{gHecke}).
\qed\\[-0.5mm]

\rem 
Analysis of (\ref{gHecke}) differs if we relax
the condition that $g(0)=0$. In this
case, by setting $\mu=\lambda$, eq.~(\ref{gHecke})
is reduced to an algebraic equation,
\be\label{gcon}
 (g(\lambda))^2 + \et \, g(\lambda) + 1 = 0 \,,
\ee
which implies that $g(\lambda)$ is a constant
function. The possible values of the constant,
i.e., the roots of (\ref{gcon}), are
the values of~(\ref{gsol}) in the limit
$\gamma\lambda \to \pm\infty$.

Now, we will apply Lemma~\ref{Hecke} in order to
obtain some information about the spectral 
decomposition of a regular $U_q(sl_2)$--invariant R--matrix.

\begin{prop}\label{r2ss}
Let $R(\lambda)$ be a $U_q(sl_2)$--invariant solution 
of the YB equation (\ref{YBE}) on $V_s^{\otimes 3}$
satisfying (\ref{gauge}) and (\ref{ans2}).
Then the second highest coefficient in its spectral 
decomposition  is given~by 
\be\label{r2s}
   r_{2s-1} (\lambda)
 = \left\{ \begin{array}{ll}
 \frac{1+\gamma\lambda}{1-\gamma\lambda}  & 
  \ \ {\rm if}\ \  q = 1 \,; \\ [1.5mm]
 \frac{[2s + \gamma\lambda]_q}%
  {[2s - \gamma\lambda]_q} &
  \ \ {\rm if}\ \  q \neq 1 \,,
 \end{array} \right. 
\ee
where $\gamma$ is an arbitrary finite constant.
\end{prop}

\proof Let $\tilde{W}_1$ denote the subspace of  
$V_s^{\otimes 3}$ which is the linear span of the
vectors
\be\label{tW}
 |1\rangle_{\one\two\three} =  |s \- 1\rangle_{\one} \, 
   |s\rangle_{\two} \, |s\rangle_{\three} \,,\quad
 |2\rangle_{\one\two\three} =  |s \rangle_{\one} \, 
   |s \- 1 \rangle_{\two} \, |s\rangle_{\three} \,,\quad
 |3\rangle_{\one\two\three} =  |s \rangle_{\one} \, 
   |s\rangle_{\two} \, |s \- 1\rangle_{\three} \,.
\ee 
{}From (\ref{YSS}) and the Clebsch--Gordan (CG) decomposition 
of~$V_s^{\otimes 2}$
(see~\cite{CG} for an explicit form of the CG coefficients),
\be
\label{ab1}
 && |2s, 2s-1  \rangle = \alpha_s \, |s\rangle |s-1\rangle 
  + \beta_s \, |s-1\rangle |s\rangle \,, \\ [1mm]
\label{ab2}
 && |2s-1, 2s-1  \rangle = \beta_s \, |s\rangle |s-1\rangle 
  - \alpha_s \, |s-1\rangle |s\rangle \,, \\ [0.5mm]
\label{albe}
 && \alpha_s = q^{s} \, (q^{2s} + q^{-2s})^{-\frac{1}{2}} 
 \,, \qquad 
 \beta_s  = q^{-s} \, (q^{2s} + q^{-2s})^{-\frac{1}{2}} \,,
\ee
we infer that $\tilde{W}_1$ is an invariant subspace of 
the YB operator for the R--matrix under consideration. 
Notice that the restrictions of $P_l^j$, 
$l={\scriptstyle \{12\}}, {\scriptstyle \{23\}}$ 
onto~$\tilde{W}_1$ vanish  if $j < 2s\- 1$. Thus,  
\be\label{RWt}
 R_l(\lambda) \Bigm|_{\tilde{W}_1} = 
 P_l^{2s} +  r_{2s\- 1} (\lambda) P_l^{2s\- 1} \,.
\ee
Furthermore, taking into account (\ref{ab1})--(\ref{ab2}) 
and introducing 
$\tilde{g}(\lambda)= 
  r_{2s\-1}(\lambda) -1$,
we observe that (\ref{RWt}) can be rewritten in 
the following form
\be\label{RWtt}
 R_l(\lambda) \Bigm|_{\tilde{W}_1} =  
 {\mathbb E} + \tilde{g}(\lambda) \, \pi_l  \,,
\ee
where 
$\pi_l$ are projectors, $\pi_l^2=\pi_l$,
given in the basis (\ref{tW})~by
\be
 \pi_{\one\two} &=& \alpha^2_s \, | 1 \rangle \, \langle 1 |
 - \alpha_s \, \beta_s \, | 1 \rangle \, \langle 2 |
 - \alpha_s \, \beta_s \, | 2 \rangle \, \langle 1 |
 + \beta^2_s \, | 2 \rangle \,  \langle 2 |    \,, \\ [0.5mm]
 \pi_{\two\three} &=& \alpha^2_s \, | 2 \rangle \, \langle 2 |
 - \alpha_s \, \beta_s \, | 2 \rangle \, \langle 3 |
 - \alpha_s \, \beta_s \, | 3 \rangle \, \langle 2 |
 + \beta^2_s \, | 3 \rangle \,  \langle 3 |  \,.
\ee
Now, noticing that
\be
 \pi_{\one\two} \, \pi_{\two\three} \, \pi_{\one\two} = 
  ( \alpha_s \, \beta_s )^2 \, \pi_{\one\two} \,, \qquad
 \pi_{\two\three} \, \pi_{\one\two} \, \pi_{\two\three} = 
  ( \alpha_s \, \beta_s )^2 \, \pi_{\two\three} \,,
\ee
we see that (\ref{RWtt}) fulfill the conditions of 
Lemma~\ref{Hecke} upon identification
$U_l= (\alpha_s \beta_s)^{-1} \pi_l$,
$g(\lambda) = \alpha_s \beta_s \, \tilde{g}(\lambda)$, and
$\et= (\alpha_s  \beta_s)^{-1}$. Substituting
this value of $\et$ into (\ref{gsol}) and recalling that
$\tilde{g}(\lambda)= r_{2s\-1}(\lambda) -1$,
we obtain~(\ref{r2s}), where we replaced $e^{\gamma\lambda}$
with $q^{2\gamma\lambda}$ for the sake of convenience
of comparison with the $q=1$ limit.
 Since we require that $R(\lambda)$
be regular, the constant $\gamma$ must be finite.
\qed

\subsection{Invariant subspaces}
The proof of Proposition~\ref{r2ss} demonstrates
that reduction of the YB operator to some invariant
subspace facilitates finding the coefficients
$r_j(\lambda)$ of an R--matrix. In what follows we will
develop this approach further exploiting available
knowledge about the CG decomposition of tensor products of
representations of~$U_q(sl_2)$. On this way, we will derive systems
of coupled functional equations similar to eq.~(\ref{gHecke})
and show that the corresponding necessary conditions
are provided by a set of coupled algebraic equations.

In Proposition~\ref{r2ss} we used that the YB 
operator~(\ref{YBO}) commutes with~$S^z_{\one\two\three}$.
Now we are going to use that 
$Y(\lambda,\mu)$ commutes with~$S^\pm_{\one\two\three}$ as well.

Let $\lfloor t \rfloor$ denote the entire part of~$t$.
Let us define the subspace $W^{(s)}_n \subset V_s^{\otimes 3}$
for \hbox{$n=0,1,\ldots,\lfloor 3s \rfloor$} as the linear 
span of highest weight vectors of spin $(3s \- n)$, i.e.,
\be\label{Wn}
  W^{(s)}_n = \{ \, \psi \in V_s^{\otimes 3} \quad \bigm| \quad
  S^+_{\one\two\three} \psi =0 \,, \quad 
  S^z_{\one\two\three} \psi = (3s \- n) \psi \, \} \,.
\ee 
Consider the following two orthonormal bases in $W^{(s)}_n$ 
(here and below 
$\FCG{\cdot}{\cdot}{\cdot}{\cdot}{\cdot}{\cdot}_q$ and 
$\FRW{\cdot}{\cdot}{\cdot}{\cdot}{\cdot}{\cdot}_q$ 
stand, respectively, for the CG coefficients and 6--$j$ symbols 
of $U_q(sl_2)$)
\be
\label{ek}
 |n;k\rangle_{\one\two\three} &=& \sum_m |m\rangle_{\one} \,
  |2s \- k , 3s \- n \- m \rangle_{\two\three} \,
 \FCG{s}{m}{2s \- k}{3s \- n \- m}{3s \- n}{3s \- n}_q
\,,\\ [1mm]
\label{ek'}
 |n;k\rangle^{\prime}_{\one\two\three} &=& \sum_m  
  |2s \- k , 3s \- n \- m \rangle_{\one\two} \, |m\rangle_{\three} \,
 \FCG{2s \- k}{3s \- n \- m}{s}{m}{3s \- n}{3s \- n}_q \,.
\ee
The basis vectors of $W^{(s)}_n$ are enumerated by 
integer $k\in I^{(s)}_n$, where  
\be\label{krange}
 I^{(s)}_n = \biggl\{ 
 \begin{array}{ll}
  0 \leq k \leq n & {\rm for} \ \  
    0 \leq n \leq 2s; \\ [0.5mm]
  n \- 2s \leq k \leq 4s \- n & {\rm for}\ \    
      2s \leq n \leq \lfloor 3s \rfloor \,.
 \end{array}  \biggm.
\ee
The sum in (\ref{ek})--(\ref{ek'}) runs over those $m$
for which the CG coefficients on the 
r.h.s.~of (\ref{ek})--(\ref{ek'}) do not vanish, i.e.,
\hbox{$(s\- n\+ k) \leq m \leq \min \bigl(s,5s\- n\- k\bigr)$}.

Consider the transition matrix, $A^{(s,n)}$, {}from the basis
(\ref{ek}) to the basis (\ref{ek'}), i.e., the orthogonal 
matrix with entries being the following scalar products 
\be
\label{defAn}
 A^{(s,n)}_{k k^{\prime} } =
  \langle n; k | n; k^{\prime} \rangle^{\prime} \,.
\ee
The transition matrix is $q$--dependent but
for compactness of notations we will not 
write the  argument $q$ explicitly unless required
by the context.

\begin{prop}\label{Aprop}

$\phantom x$\par\noindent
i) entries of $A^{(s,n)}$ are expressed in terms of 6--$j$ symbols
of $U_q(sl_2)$ as follows
\be\label{A1}
 A^{(s,n)}_{k k^{\prime} }  = 
 (-1)^{2s-n} \sqrt{[4s\- 2k \+ 1]_q [4s \- 2k' \+ 1]_q} \
 \FRW{s}{s}{s}{3s \- n}{2s \- k}{2s \- k'}_q \,.
\ee
ii) $A^{(s,n)}$ is self--dual in $q$, 
\be
\label{A2}
 A^{(s,n)}_q = A^{(s,n)}_{q^{-1}} \,.
\ee
iii) $A^{(s,n)}$ is orthogonal, symmetric, and
 coincides with its inverse ($t$ denotes the matrix 
transposition),
\be
\label{A3}
 A^{(s,n)} =  \bigl( A^{(s,n)} \bigr)^t = 
  \bigl( A^{(s,n)} \bigr)^{-1} \,.
\ee
As a consequence, the only eigenvalues of 
$A^{(s,n)}$ are~$\pm 1$.
\par\noindent
iv) Transition matrices enjoy the following ``spin--level duality''
relations
\be
\label{A6a}
  A^{(s,1)}_q  &=& A^{(\frac{1}{2},1)}_{q^{2s}}  \,, \\ [1mm]
\label{A6b}
  A^{(s,n)}_q  &=& 
  A^{(2s-\frac{n}{2}, 6s-2n)}_{q}  \,, 
\ee
where $n \leq 2s$.
\end{prop}
Explicit formulae for entries of matrix $A^{(s,n)}$ and a proof 
of its properties listed above are given in Appendices~B and~C.

\rem
It follows from $iii$) and $iv$) that 
$\frac{1}{2} \bigl({\mathbb E} \pm A^{(s,n)} \bigr)$
are projectors of ranks~$n_\pm$. 
In particular, for $n\leq 2s$, we have 
$ n_+ = \lfloor \fr{n}{2} \+ 1 \rfloor$,
$ n_- = \lfloor \fr{n+1}{2} \rfloor $.

Properties of the transition matrix given in  
Proposition~\ref{Aprop} make it an efficient tool for
dealing with restrictions of $U_q(sl_2)$--invariant
operators to subspaces~$W^{(s)}_n$. As a simple example,
let us prove the following well--known statement.
\begin{lem}\label{PPP}
The following identity holds on~$V_s^{\otimes 3}$
\be\label{ppp1}
  P_{\two\three}^0 \, P^j_{\one\two} \, P^0_{\two\three} = 
 \frac{[2j\+1]_q}{[2s\+1]_q^2} \, P^0_{\two\three} \,.
\ee
\end{lem}

\proof  Observe that  
$P^0_{\one\two} \bigl|_{W^{(s)}_{n}}$ and 
$P^0_{\two\three} \bigl|_{W^{(s)}_{n}}$
vanish for all $n$ except $n=2s$. Thus it suffices 
to prove (\ref{ppp1}) when it is restricted
to~$W^{(s)}_{2s}$. Denote 
$p^j= P^j_{\two\three} \bigl|_{W^{(s)}_{2s}}$.
In the basis (\ref{ek}) we have 
$p^j_{ab}=\delta_{a,2s-j}\delta_{b,2s-j}$. Therefore,
in this basis, the l.h.s.~of~(\ref{ppp1})
acquires the following form
\be\label{ppp2}
 p^0 \, A^{(s,2s)} \, p^j \, A^{(s,2s)} \, p^0
 = \Bigl( A^{(s,2s)}_{2s\-j,2s} \Bigr)^2 \, p^0 \,.
\ee
The value of $A^{(s,2s)}_{2s\-j,2s}$ is easily 
computable (see eq.~(\ref{Akn}) in Appendix~B)
and its square yields the scalar coefficient on the 
r.h.s.~of~(\ref{ppp1}).

\subsection{Reduced Yang--Baxter equations}
 
Eqs.~(\ref{YSS}) imply that $W^{(s)}_n$ is an invariant subspace 
for the YB operator~(\ref{YBO}). Let us introduce the
{\em reduced} YB operator:
$Y_n(\lambda,\mu)=Y(\lambda,\mu)\!\bigm|_{W^{(s)}_n}$ 
(the restriction of $Y(\lambda,\mu)$ onto~$W^{(s)}_n$).
Notice that restrictions of $P_l^j$ to $W^{(s)}_n$
are diagonal in the bases (\ref{ek'})
and (\ref{ek}) for $l={\scriptstyle \{12\}}$ and 
$l={\scriptstyle \{23\}}$, 
respectively. Moreover, they vanish unless 
\be\label{jWn}
 |2s - n| \leq j \leq \min(2s,4s-n) \,.
\ee
Therefore, in the basis (\ref{ek'}), 
$R_l(\lambda)\!\bigm|_{W^{(s)}_n}$ are represented as
\be
 & R_{\one\two}(\lambda) \Bigm|_{W^{(s)}_n} =
 A^{(s,n)} \, D(\lambda) \, \bigl( A^{(s,n)} \bigr)^{-1}
\,, \qquad
 R_{\two\three}(\lambda) \Bigm|_{W^{(s)}_n} = D(\lambda) \,, \\
 & D_{k k'}(\lambda) = \delta_{k k'} \, r_{2s-k}(\lambda) \,,&
\ee
where $k \in I_n^{(s)}$ as specified in~(\ref{krange}).
Whence, taking the property (\ref{A3}) into account,
we conclude that $Y_n(\lambda,\mu)$ acquires
the following form in the basis (\ref{ek'}) 
\be
 \label{Yn}
  Y_n(\lambda,\mu) &= &
 A^{(s,n)} \, D(\lambda-\mu) \, A^{(s,n)} \, 
 D(\lambda) \, A^{(s,n)} \, D(\mu) \\ [0.5mm]
\nonumber && 
 -\, D(\mu) \, A^{(s,n)} \, D(\lambda) \,
 A^{(s,n)} \, D(\lambda-\mu) \, A^{(s,n)}  .
\ee
The corresponding {\em reduced} YB equation reads
\be\label{RYB}
 A^{(s,n)} \, D(\lambda-\mu) \, A^{(s,n)} \, 
 D(\lambda) \, A^{(s,n)} \, D(\mu) =
 D(\mu) \, A^{(s,n)} \, D(\lambda) \,
 A^{(s,n)} \, D(\lambda-\mu) \, A^{(s,n)} \,.
\ee
It is the condition of vanishing of the YB operator
(\ref{YBO}) on~$W^{(s)}_n$. Observe that  (\ref{A3}) implies 
that the reduced YB operator is antisymmetric,
$ \bigl( Y_n(\lambda,\mu) \bigr)^t = - Y_n(\lambda,\mu)$.
Therefore independent relations contained in (\ref{RYB}) are
\be
\nonumber
 && \sum_{i,j \in  I^{(s)}_n} 
 r_{2s-i} (\lambda \- \mu) \, r_{2s-j} (\lambda) \,
 A^{(s,n)}_{ij} \,
 \Bigl( r_{2s-a} (\mu) \, A^{(s,n)}_{aj} \,
 A^{(s,n)}_{ib}   \\
\label{rYB}
 && \qquad\qquad\qquad
 - r_{2s-b} (\mu) \, A^{(s,n)}_{ai} \,
 A^{(s,n)}_{jb} \Bigr)= 0 \,, \qquad
 a < b \,, \quad a,b \in  I^{(s)}_n \,.
\ee

Let us emphasize that eqs.~(\ref{rYB}) at the level~$n$ 
ensure, thanks to commutativity of the YB operator with 
$S^-_{\one\two\three}$, that the YB operator vanishes 
not only on the subspace $W^{(s)}_n$ but also on the 
larger subspace that is spanned by {\em all} vectors
obtained by acting on~$W^{(s)}_n$ with 
$(S^-_{\one\two\three})^m$, $m=0,\ldots,6s-2n$.
(This picture  resembles closely the structure of
eigenvectors in the algebraic Bethe ansatz,
see~\cite{Fad1} for a review).
Thus, the set of reduced YB equations (\ref{rYB}),
$n=1,\ldots,\lfloor 3s \rfloor$, is less
overdetermined than the initial YB equation (\ref{YBE}) 
containing ${\rm dim}\, V_s^{\otimes 3}=
(2s\+1)^3$ functional equations (although some of them 
are in general not independent). 
However, even this set of equations is still overdetermined. 
Indeed, (\ref{rYB}) at level $n\leq 2s$ involves $r_j(\lambda)$
with $j=2s\-n,\ldots,2s$. Therefore, it suffices to solve 
(\ref{rYB}) for $n=1,\ldots,2s$ to determine all coefficients 
$r_j(\lambda)$ of an R--matrix. But these coefficients also
have to satisfy the remaining reduced YB equations for 
$n=2s\+1,\ldots,\lfloor 3s \rfloor$.

\rem For $s = \frac 12$ we have $2s\=\lfloor 3s \rfloor \=1$
and therefore the corresponding set of reduced YB equations
is not overdetermined. Indeed, in this case (\ref{rYB})
contains only one independent equation. A slightly less 
trivial remark is that the set of reduced YB equations is 
not overdetermined for $s=1$ as well (see the proof of
Proposition~\ref{Rs1}). 

Now, as an immediate application of the reduced YB equation 
technique, let us prove the following statement. 
\begin{prop}\label{Hc}
Let $R(\lambda)$ be a $U_q(sl_2)$--invariant solution 
of the YB equation (\ref{YBE}) on $V_s^{\otimes 3}$
for a spin~$s\geq \frac{n}{2}$, $n\in {\mathbb Z}_+$
satisfying (\ref{gauge}) and (\ref{ans2}). 
Suppose that $n$ highest coefficients
in its spectral decomposition coincide, 
$r_{2s}(\lambda)=r_{2s\-1}(\lambda)=
 \ldots=r_{2s\-n\+1}(\lambda)=1$.
Then 
\be\label{r3a}
 r_{2s\-n} (\lambda) = 1 + \et \, g(\lambda) \,,
\ee
where $g(\lambda)$ is given by (\ref{gsol}) with
\be\label{d0s}
 \et = \frac{[2s\-n]!}{[2s]!} \, 
  \frac{[4s\-n\+1]!}{[4s\-2n\+1]!} \,.
\ee
Here the $q$--factorial is defined as  
$[n]! =\prod_{k=1}^n [k]_q$ for $n\in {\mathbb Z}_+$ 
and~$[0]!=1$.
\end{prop} 

\proof The corresponding reduced YB equation
(with $n$ in (\ref{RYB}) being the same as in~(\ref{r3a}))
multiplied {}from the left by $A^{(s,n)}$ can
be regarded as the YB equation (\ref{YBE}) for 
$R_{\one\two}(\lambda)=D(\lambda)$ and
$R_{\two\three}(\lambda)=
 A^{(s,n)} \, D(\lambda) \, A^{(s,n)}$.
Further, we notice that 
\be
 D(\lambda) = {\mathbb E} + \tilde{g}(\lambda) \, \pi
 \,, \qquad
 A^{(s,n)} \, D(\lambda) \, A^{(s,n)} =
 {\mathbb E} + \tilde{g}(\lambda) \, \pi^{\prime} \,,
\ee
where $\tilde{g}(\lambda)=r_{2s\-n} (\lambda)-1$,
$\pi$ is a matrix such that $\pi_{ab}=\delta_{an} \delta_{bn}$,
$a,b=0,\ldots,n$,
and \hbox{$\pi^{\prime}=A^{(s,n)} \, \pi \, A^{(s,n)}$}.
Obviously, $\pi$ and $\pi^{\prime}$ are projectors
of rank one. Moreover, a computation similar to~(\ref{ppp2})
shows that
\be
 \pi \, \pi^{\prime} \, \pi = \et^{-2} \, \pi \,, \qquad
 \pi^{\prime} \, \pi \, \pi^{\prime} = 
   \et^{-2} \, \pi^{\prime} \,, \qquad
 \et = \bigl|A^{(s,n)}_{nn}\bigr|^{-1} \,.
\ee
Hence (\ref{r3a}) follows by invoking Lemma~\ref{Hecke}
upon the identification
$U_{\one\two}=\et \, \pi$,
$U_{\two\three}=\et \, \pi^{\prime}$, and
$\tilde{g}(\lambda) = \et \, g(\lambda)$.
Explicit form of $\et$ given in (\ref{d0s}) 
is easily obtained~from~(\ref{Akn}).
\qed\\[-0.5mm]

This Proposition generalizes both Lemma~\ref{Hecke} and 
Proposition~\ref{r2ss}.
For $n=1$, eq.~(\ref{d0s}) yields $\et=q^{2s}\+q^{-2s}$
and we recover the case of Proposition~\ref{r2ss}.
For $n=2s$, eq.~(\ref{d0s}) yields $\et=[2s\+1]_q$;
the corresponding R--matrix is given by~(\ref{qTL})
which is a particular example covered by Lemma~\ref{Hecke}
(cf.~Remark~3).
It is not clear whether an example of R--matrix
with such coefficients $r_{j}(\lambda)$ as described in 
Proposition~\ref{Hc} exists for $n\neq 1$ and $n\neq 2s$.
Nevertheless this proposition is useful for the analysis 
of solutions of the YB equation (see the next section).
Another statement useful for this analysis reads as follows.
\begin{prop}\label{rhalf}
Let $R(\lambda)$ be a $U_q(sl_2)$--invariant solution 
of the YB equation (\ref{YBE}) on $V_s^{\otimes 3}$
for a half--integer spin~$s\geq \fr{3}{2}$
satisfying (\ref{gauge}) and~(\ref{ans2}).
Then the coefficients $r_{s\-\frac{1}{2}}(\lambda)$ and
$r_{s+\frac{1}{2}}(\lambda)$ in its spectral 
decomposition are related as follows 
\be\label{rshalf}
 \frac{r_{s-\frac{1}{2}} (\lambda)}{r_{s+\frac{1}{2}} (\lambda)} 
 = \left\{ \begin{array}{ll}
 \frac{1+\gamma\lambda}{1-\gamma\lambda}  & 
  \ \ {\rm if}\ \  q = 1 \,; \\ [1.5mm]
 \frac{[s\+ \frac{1}{2} + \gamma\lambda]_q }%
  {[s\+ \frac{1}{2} - \gamma\lambda]_q }  &
  \ \ {\rm if}\ \  q \neq 1 \,,
 \end{array} \right. 
\ee
where $\gamma$ is an arbitrary finite constant.
\end{prop}

\proof Matrix $D(\lambda)$ in the reduced YB equation 
can be multiplied by an arbitrary function
$\varphi(\lambda)$ analytic in a neighbourhood of
$\lambda=0$ and satisfying 
$\varphi(\lambda)\varphi(-\lambda)=1$.
Therefore, in (\ref{RYB}) for $n=3s\- \frac{1}{2}$,
we can choose $D(\lambda)={\rm diag}\, (1,g(\lambda))$,
where
$g(\lambda)=\frac{r_{s-\frac{1}{2}} (\lambda)}%
{r_{s+\frac{1}{2}} (\lambda)}$. Further, by
the duality relation~(\ref{A6b}), we have 
$ A^{(s,3s\- \frac{1}{2})}  = 
  A^{(\frac{s}{2}\+\frac{1}{4},1)}$
for half--integer spins $s\geq \fr{3}{2}$.
Whence, applying Proposition~\ref{r2ss},
we conclude that 
$g(\lambda)= r_{2s'-1}(\lambda)$, where
$s'=\frac{s}{2} + \frac{1}{4}$.
\qed

\subsection{Necessary conditions}

Differentiating (\ref{rYB}) w.r.t.~$\mu$, 
setting $\mu=\lambda$, and taking into account the regularity 
condition $D(0)={\mathbb E}$, we derive the following system 
of equations (prime denotes derivative
w.r.t.~the spectral parameter):
\be
\label{dDA}
 && \sum_{i,j \in  I^{(s)}_n} 
 r_{2s-i}^\prime (0) \, r_{2s-j} (\lambda) \, A^{(s,n)}_{ij} \,
 \Bigl( r_{2s-a} (\lambda) \, A^{(s,n)}_{aj} \,
 A^{(s,n)}_{ib} - 
 r_{2s-b} (\lambda) \, A^{(s,n)}_{ai} \,
 A^{(s,n)}_{jb} \Bigr) \\
\nonumber
&& \quad =  A^{(s,n)}_{ab} \,
 \bigl( r_{2s-a}^\prime (\lambda) \, r_{2s-b} (\lambda) -
 r_{2s-b}^\prime (\lambda) \, r_{2s-a} (\lambda) \bigr) \,,
 \qquad  a < b \,, \quad a,b \in  I^{(s)}_n \,. 
\ee
Here we have carried out the summation on the r.h.s.\ 
by using that $\bigl(A^{(s,n)} A^{(s,n)}\bigr)_{ab}
 = \delta_{ab}$.
It is important to emphasize here that although (\ref{dDA}) 
contain derivatives, they are actually linear {\em algebraic} 
equations on $r_j (\lambda)$ for $j\neq a,b$.

It is easy to check that (\ref{dDA}) is satisfied trivially
for~$\lambda=0$. Therefore, let us look at higher
order terms in the expansion of (\ref{dDA}) about 
$\lambda=0$. Denote $r_{2s-a}' (0) \equiv \xi_a$. In the 
first order in $\lambda$, the summation over $i,j$ can be 
carried out, and we obtain the conditions
\be\label{dd0}
 A^{(s,n)}_{ab} \, \bigl( r_{2s-a}'' (0) - 
  r_{2s-b}'' (0) \bigr)
 = A^{(s,n)}_{ab} \, \bigl(  \xi_a^2 -
 \xi_b^2 \bigr) \,,
\ee
that are always satisfied, because the unitarity (\ref{uni})
implies that
\be\label{dd}
 r_{2s-a}'' (0) = \xi_a^2 \,.
\ee
In the second order in $\lambda$, eqs.~(\ref{dDA})
turn into a system of algebraic equations
\be  
\label{ddd}
 && \sum_{i,j \in  I^{(s)}_n}  \xi_i \, \xi_j^2 \,
 A^{(s,n)}_{ij} \,
 \Bigl( A^{(s,n)}_{aj} \, A^{(s,n)}_{ib} - 
 A^{(s,n)}_{ai} \, A^{(s,n)}_{jb} \Bigr) \\
\nonumber
 && + \, \bigl( \xi_a - \xi_b \bigr) \,
 \sum_{i,j \in  I^{(s)}_n}   \xi_i \, \xi_j  \,
 A^{(s,n)}_{ij} \,
 \Bigl( A^{(s,n)}_{aj} \, A^{(s,n)}_{ib} + 
 A^{(s,n)}_{ai} \, A^{(s,n)}_{jb} \Bigr) \\
\nonumber
&& =  A^{(s,n)}_{ab} \,
 \Bigl( r_{2s-a}''' (0)  -  r_{2s-b}''' (0) 
 - \xi_a^3 + \xi_b^3 + \xi_a^2 \, \xi_b - \xi_b^2 \, \xi_a
 \Bigr) \,,
 \qquad  a < b \,, \quad a,b \in  I^{(s)}_n \,. 
\ee

Remarkably, equations (\ref{dDA}) and 
(\ref{ddd}) can be solved in a recursive way. 
Let us provide the corresponding algorithm.
We start with the level~$n=1$, where we have $r_{2s}(\lambda)=1$
and $r_{2s-1}(\lambda)$ is given by~(\ref{r2s}).
Now, suppose we have found $r_j(\lambda)$, $j=2s\-n,\ldots,2s$
that solve eqs.~(\ref{dDA}) for a given level~$n<2s$. 
Then (\ref{dDA}) and (\ref{ddd}) at the
level $(n\+1)$ allow us to express $r_{2s\-n\-1}(\lambda)$ 
{\em algebraically} in terms of the previously
found~$r_j(\lambda)$. Indeed, since we already
know $\xi_j$ for $j=0,\ldots,n$, eqs.~(\ref{ddd})
for $2s\-n \leq a<b \leq n $ turn into quadratic equations
w.r.t.~$\xi_{n\+ 1}$. Solving them and substituting
the found values of $\xi_{n\+ 1}$ into~(\ref{dDA})
for $2s\-n \leq a<b \leq n $, we obtain a system of linear 
equations on $r_{2s\-n\-1}(\lambda)$. Finding all possible
solutions to this system completes the $(n\+1)$--th
step of recursion. Continuing this procedure up to $n=2s$,
we will obtain all possible solutions for all~$r_j(\lambda)$
and thus construct {\em all} possible ans\"atze for regular
$U_q(sl_2)$--invariant R--matrices of spin~$s$.
Next, since (\ref{dDA}) provide necessary but not sufficient
conditions, we have to check which of these ans\"atze
indeed satisfy the YB equation (\ref{YBE}) or, alternatively,
the reduced YB equations (\ref{rYB}) for all $n$ up 
to~$\lfloor 3s \rfloor$.

\subsection{Spin chain Hamiltonians and reconstruction
of R--matrices}

The utmost importance of the YB equation in the quantum 
inverse scattering method (see~\cite{KS,Fad1} for a review) 
is due to the fact that its solutions can be used to
construct families of quantum integrals of motion 
in involution. In particular, regular solutions
of the YB equation allow to construct {\em local}
integrals of motion for lattice models. For
the R--matrix of type (\ref{ans1}), the first of 
these integrals,
\be\label{H}
 {\cal H} = \sum_k H_{k,k+1} \,, \qquad
 H = \partial_\lambda R(\lambda) 
 \Bigm|_{\lambda=0}  = 
  \sum_{j=0}^{2s-1} \xi_{2s-j} \, P^j \,,
\ee
is usually regarded as a Hamiltonian of a spin $s$ 
magnetic chain with the nearest neighbour interaction. 
Here \hbox{$H \in {\rm End}\ V_{s}^{\otimes 2}$} and
${\cal H} \in {\rm End}\ V_{s}^{\otimes L}$, where
$L$ is the number of lattice sites. Notice
that in (\ref{H}) we took into account the
normalization condition~(\ref{gauge}), which implies
$\xi_{0}=0$ (this fixes the choice of the additive
constant in the Hamiltonian). \hbox{R--matrices}
equivalent in the sense of transformation (\ref{ren})
yield Hamiltonians related simply by rescaling
$H \to \gamma H$; we will regard such
Hamiltonians as equivalent.

It is important to remark that, as it was observed
in~\cite{Ke}, for regular solutions of the YB equation  
different Hamiltonians correspond to inequivalent R--matrices.
In the present context this statement can be formulated
as follows.

\begin{lem}\label{HR}
Let $R^{(1)}(\lambda)$ and $R^{(2)}(\lambda)$ be two 
solutions of the YB equation (\ref{YBE}) on 
$V_{s}^{\otimes 3}$ satisfying 
(\ref{ans1}), (\ref{gauge}), and (\ref{ans2}). The 
corresponding Hamiltonians given by~(\ref{H})
coincide, $H^{(1)}=H^{(2)}$, if and only if
$R^{(1)}(\lambda)=R^{(2)}(\lambda)$.
\end{lem}
\proof The ``if'' part is obvious. Further,  
Theorem~3 in~\cite{Ke} asserts that
if the Hamiltonians corresponding to two
regular R--matrices analytic in a neighbourhood
of~$\lambda=0$ coincide, then 
$R^{(1)}(\lambda)=\varphi(\lambda) \,R^{(2)}(\lambda)$,
where the scalar function $\varphi(\lambda)$ is analytic 
in a neighbourhood of~$\lambda=0$ and 
satisfies the condition $\varphi(0)=1$. 
In the case under consideration, analyticity
of $r_j(\lambda)$ along with condition (\ref{gauge}) 
imply that $\varphi(\lambda)=1$.
\qed\\[-0.5mm]

\rem 
The algorithm described at the end of the previous subsection
complements this Lemma with a constructive procedure that
allows us to {\em reconstruct} the R--matrix {}from a 
given Hamiltonian. Indeed, if we know a Hamiltonian
in the form~(\ref{H}), i.e., we know all~$\xi_j$,
then we can solve (\ref{dDA}) recursively starting
with $r_{2s}(\lambda)=1$ thus recovering all 
the coefficient $r_j(\lambda)$ of the corresponding
regular $U_q(sl_2)$--invariant R--matrix. 
In contrast to the general situation, Lemma~\ref{HR}
guarantees that the resulting set of $r_j(\lambda)$ will
be unique.

\section{Analysis of reduced Yang--Baxter equations}

\subsection{Asymptotic solutions}
Let us remark that the technique described above applies
in the limit $\lambda\rightarrow\infty$ as well.
In this limit, assuming that 
$\check{R}^{\pm 1}= 
\lim_{\lambda\rightarrow\pm\infty}R(\lambda)$
exist, the YB equation~(\ref{YBE}) turns into 
\be\label{YBc}
 \check{R}_{\one\two} \, \check{R}_{\two\three} \, 
 \check{R}_{\one\two} = \check{R}_{\two\three} \, 
 \check{R}_{\one\two} \, \check{R}_{\two\three} \,.
\ee 
Denote $d_j = \lim_{\lambda\rightarrow +\infty}r_j(\lambda)$,
so that $\check{R} = \sum_{j=0}^{2s} d_j P^j$.
Condition (\ref{gauge}) implies that~$d_{2s}\= 1$.
Taking the limit $\lambda\rightarrow\infty$ 
in (\ref{RYB}), we obtain a set of algebraic equations
on the coefficients~$d_j$,
\be\label{Ync}
 A^{(s,n)} \, D \, A^{(s,n)} \, 
 D \, A^{(s,n)} \, D =
 D \, A^{(s,n)} \, D \,
 A^{(s,n)} \, D \, A^{(s,n)} \,.
\ee
Here $n=1,\ldots,\lfloor 3s \rfloor$, and
$D_{k k'}= \delta_{k k'} \, d_{2s-k}$, 
where $k \in I_n^{(s)}$.
Analogously to the spectral dependent case,
independent equations contained in (\ref{Ync}) are
\be\label{ddAAA}
 (d_a - d_b ) \, \sum_{i,j \in  I^{(s)}_n}  d_i \, d_j \,
 A^{(s,n)}_{ij} \, A^{(s,n)}_{ai} \, A^{(s,n)}_{jb} = 0 \,,
 \qquad  a < b \,, \quad a,b \in  I^{(s)}_n \,.
\ee
This system of equations can be solved in a recursive way
with the help of the algorithm described at the 
end of subsection~2.4. In this context, it is worth 
noticing that one particular solution is known
a--priory, namely 
\be
\label{dj}
 d_j = (-1)^{2s-j} \, q^{2s(2s+1)-j(j+1)} \,,
\ee
which corresponds to (\ref{serB}) in the limit
$q^\lambda\rightarrow\infty$.

\subsection{Reduced YB for $n=1$ and $n=2$}

The explicit form of $n=1$ transition matrix is
\be
\label{n1}
 A^{(s,1)} =  \frac{1}{ \{2s\}_q } \, 
 \biggl( \begin{array}{cc} 1 & \sqrt{1\+ \{4s\}_q } \\
 \sqrt{1\+ \{4s\}_q } & -1  \end{array} \biggr) \,,
\ee
where we denoted $\{t\}_q = q^t + q^{-t}$. In this case 
eq.~(\ref{dDA}) contains only one equation.  Taking into 
account that $r_{2s}(\lambda)=1$ and making substitution
$r_{2s-1}(\lambda) = 1 + \{2s\}_q \, g(\lambda) $,
it is easy to see that this equation coincides
with~(\ref{gdif}). Hence we recover the same expression 
for $r_{2s-1}(\lambda)$ as in Proposition~\ref{r2ss}.
Accordingly, eq.~(\ref{ddAAA}) is either satisfied 
trivially if $d_{2s-1}=1$ or it represents a quadratic
equation with roots $d_{2s-1}=-q^{\pm 4s}$.

The explicit form of $n=2$ transition matrix is
\be
\label{n2}
 A^{(s,2)} = 
 \left( \begin{array}{ccc} 
 \frac{[2s\-1]_q}{  \{2s\}_q \, [4s\-1]_q} &
 \rho_s 
 \sqrt{ \frac{[2]_q\, [6s\-1]_q}{ 
 \{2s\}_q \, [4s\-1]_q} }  & 
 \rho_s 
 \frac{ \sqrt{  [6s\-1]_q \, [6s\-2]_q }}{ 
 [4s\-1]_q  } \\ [2.5mm]
 \rho_s 
 \sqrt{ \frac{[2]_q\, [6s\-1]_q}{ 
 \{2s\}_q \, [4s\-1]_q} } & 
 \{ 4s\-1 \}_q \, (\rho_s)^2 &
 - \rho_s  
  \sqrt{ \frac{ [2]_q \, [6s\-2]_q}{ 
 \{2s\-1\}_q \, [4s\-1]_q  }} \\ [2.5mm]
 \rho_s 
 \frac{ \sqrt{  [6s\-1]_q \, [6s\-2]_q }}{ 
 [4s\-1]_q  } &
 - \rho_s  
  \sqrt{ \frac{ [2]_q \, [6s\-2]_q}{ 
 \{2s\-1\}_q \, [4s\-1]_q  }} &
 \frac{ [2s]_q}{ \{2s \- 1\}_q [4s \- 1]_q} 
 \end{array} \right) , 
\ee
where we introduced
$\rho_s = \bigl( \{2s\-1\}_q \, \{2s\}_q \bigr)^{-\frac{1}{2}}$.
For computations, it is useful to observe the 
following identities:
\be 
 A_{01}^{(s,2)} \, A_{12}^{(s,2)} =
 A_{02}^{(s,2)} \, \bigl( A_{11}^{(s,2)} -1 \bigr) \,,\qquad
 A_{11}^{(s,2)} = 1 \- [2]_q \, (\rho_s)^2 \,. 
\ee

Analysis of (\ref{rYB}) splits into two cases:
$r_{2s-1}(\lambda)=1$ and $r_{2s-1}(\lambda) \neq 1$. The
former case is covered by Proposition~\ref{Hc}, which yields
\be\label{tl2}
 r_{2s-2}(\lambda) =
 \frac{b^2 e^{\lambda} -1}{b^2 - e^{\lambda}} \,, \qquad 
 b+b^{-1}=\frac{[4s\-1]_q \{2s\-1\}_q }{[2s]_q} \,.
\ee
In the latter case, without loss of generality we can
choose $\gamma=1$ in (\ref{r2s}), which corresponds to
\be\label{x1}
 r_{2s-1}(\lambda)=\frac{[2s\+\lambda]_q }{[2s \- \lambda]_q }
 \,, \qquad 
 d_{2s-1}=-q^{4s} \,, \qquad
 \xi_{1}= \varkappa_q \, \frac{ \{2s\}_q }{[2s]_q}
 \,, \qquad
 \varkappa_q \equiv  \frac{2\log q}{q\-q^{-1}} \,.
\ee
In this case, one finds easily that
the three equations contained in (\ref{ddAAA}) 
for $n=2$ have only one common root, namely
\be
\label{d2s3}
 d_{2s-2} = q^{8s-2}  \,.
\ee
Substituting $\xi_{0}=0$ and $\xi_{1}$ given by (\ref{x1})
into (\ref{ddd}), we obtain a quadratic equation 
on~$\xi_{2}$ with roots given by
\be\label{sd3}
 \xi_{2} = \frac{2 \varkappa_q \,[4s\-1]_q}%
 {[2s\-1]_q\,[2s]_q} \,, \qquad
  \xi_{2} = \varkappa_q \, (q\-q^{-1})^2 \, 
 \frac{[4s-1]_q}{ \{4s-1\}_q } \,.
\ee
The corresponding solutions of (\ref{dDA}) are given by
\be\label{rn2}
 r_{2s-2}(\lambda) = 
 \frac{[2s\+\lambda]_q }{[2s \- \lambda]_q } 
 \frac{[2s\- 1\+\lambda]_q }{[2s\-1 \- \lambda]_q }
 \,, \qquad
 r_{2s-2}(\lambda) = 
 \frac{ \{4s\-1 \+\lambda\}_q }{\{4s\-1\-\lambda\}_q  } \,.
\ee
It is straightforward to verify that both these solutions 
satisfy the $n=2$ level reduced YB equation~(\ref{rYB}).

\begin{prop}\label{Rs1}
For a generic $q$ and spin $s=1$, inequivalent regular 
$U_q(sl_2)$--invariant solutions of the YB equation 
(\ref{YBE}) satisfying the condition~(\ref{gauge})
are exhausted by the following three types
\be
\label{s1a}
&& R(\lambda) = P^2 + P^1 +
 \frac{b^2 e^{\lambda} -1}{b^2 - e^{\lambda}} P^0 
 \,, \quad b+b^{-1}=[3]_q \,, \\ [0.5mm]
\label{s1b}
&&  R(\lambda) = P^2 + 
 \frac{[2 \+ \lambda]_q}{[2 \- \lambda]_q } \, P^1 
 + \frac{[2 \+ \lambda]_q}{[2 \- \lambda]_q } 
 \frac{[1 \+ \lambda]_q} {[1 \- \lambda]_q } \, P^0 \,, \\ [0.5mm]
\label{s1c}
&& R(\lambda) = P^2 + 
 \frac{[2 \+ \lambda]_q}{[2 \- \lambda]_q } \, P^1 +
 \frac{ \{3 \+\lambda\}_q }{ \{ 3 \-\lambda\}_q  } \, P^0 \,.
\ee
\end{prop}

\proof For $n \leq 2s$ we have ${\rm dim}\, W^{(s)}_n=(n\+ 1)$.
As it has already been mentioned, reduced YB equation~(\ref{RYB})  
ensures vanishing of the YB operator on all vectors of the form 
$(S_{\one\two\three}^-)^m W^{(s)}_n$, i.e., on a subspace of
dimension 
$\Delta^{(s)}_n = (6s\-2n\+1)\,{\rm dim}\, W^{(s)}_n$.
In particular, we have 
$\Delta^{(1)}_0 \+ \Delta^{(1)}_1 \+ \Delta^{(1)}_2 = 26$, 
which means that the level $n=3$ reduced YB equation
is satisfied automatically, since the corresponding 
subspace $W^{(1)}_3$ is one dimensional. Therefore, 
for $s=1$, the $n=1,2$ reduced YB equations provide not 
only necessary but also sufficient conditions.
As we have shown above in this subsection, solutions to
these equations are exhausted by (\ref{tl2}) and~(\ref{rn2})
which for $s=1$ yield (\ref{s1a})--(\ref{s1c}). 
\qed \\[-0.5mm]

Thus, for spin $s=1$, all three inequivalent $sl_2$--invariant 
R--matrices (\ref{ser0d}), (\ref{ser0b}), and (\ref{ser0a})
have $U_q(sl_2)$--invariant counterparts. Two of them 
belong to the well--known types (\ref{serB}) and (\ref{qTL}).
The last one, (\ref{s1c}), appears to be rather an exceptional
case; it was found previously~\cite{Jo} by means of 
Baxterization of the Birman--Wenzl--Murakami algebra.

\subsection{Reduced YB for $n=3$}

For $n=3$ the transition matrix has 12 entries
given by (\ref{Akn})--(\ref{A0k}) and the remaining
four entries are 
\be\label{An3a}
 && A^{(s,3)}_{11} = 
 \frac{ [2]_q [2s\-1]_q [6s\-2]_q - ([2s\-2]_q)^2}%
 { \{2s\-1\}_q [4s\-3]_q [4s]_q } \,,  \\
\label{An3b}
 && A^{(s,3)}_{12} = A^{(s,3)}_{21} =
  \frac{ [2s\-2]_q }{ [4s\-2]_q }
 \bigl( [6s\-2]_q - [2]_q [2s\-1]_q \bigr)
 \sqrt{\frac{ [2s\-1]_q [6s\-3]_q}
 {[4s\-4]_q [4s\-3]_q [4s\-1]_q [4s]_q} } ,
 \hphantom{XXX} \\
\label{An3c}
 && A^{(s,3)}_{22} = 
 \frac{ [2s\-2]_q - [2]_q [6s\-3]_q }%
 { \{2s\-2\}_q \{2s\-1\}_q  [4s\-1]_q  } \,.
\ee

\begin{prop}\label{Rs32}
For a generic $q$ and spin $s=\frac{3}{2}$, 
inequivalent regular $U_q(sl_2)$--invariant solutions 
of the YB equation (\ref{YBE}) satisfying the
condition~(\ref{gauge}) are exhausted by the 
two types given by (\ref{serB}) and (\ref{qTL}).
\end{prop}

\proof Let us analyse the spectral decomposition of possible 
spin $\frac{3}{2}$ solutions to the reduced YB equations for 
$n=1,2,3$. The first possibility is 
$r_2(\lambda)=r_1(\lambda)=1$, in which case $r_0(\lambda)$ 
is determined by Proposition~\ref{Hc}; the corresponding
R--matrix is given by~(\ref{qTL}). Next, the case
$r_2(\lambda)=1$, $r_1(\lambda) \neq 1$
is covered by the same proposition for $n=2$ and 
$r_1(\lambda)$ is given by~(\ref{tl2}). 
However, this case is ruled out, because (\ref{tl2}) 
for $s=1$ is incompatible with the statement of 
Proposition~\ref{rhalf} which requires~$b=q^2$. 
In the remaining case, $r_2(\lambda) \neq 1$,  
without loss of generality we can choose $\gamma=1$ in 
(\ref{r2s}), which yields 
$r_2(\lambda) =\frac{[3 \+ \lambda]_q}{[3 \- \lambda]_q }$,
and, according to the analysis carried out in the
previous subsection, $r_1(\lambda)$ is given by 
one of the expressions in~(\ref{rn2}). However, the second 
form in~(\ref{rn2}) is ruled out again as incompatible 
with Proposition~\ref{rhalf}. Thus, we are left with 
\be\label{x12}
 r_2(\lambda)=\frac{[3 \+ \lambda]_q}{[3 \- \lambda]_q }
 \frac{[2 \+ \lambda]_q}{[2 \- \lambda]_q} \,, \qquad
 \xi_0 = 0 \,, \qquad
 \xi_1 = \varkappa_q \, \frac{ \{3\}_q }{[3]_q}
 \,, \qquad 
 \xi_{2}= 2\varkappa_q \, \frac{ [5]_q }{[2]_q [3]_q}
\,.
\ee
Substituting these values into (\ref{ddd}) for $n=3$ 
and $s=\frac{3}{2}$, we obtain a system of three quadratic 
equations on~$\xi_{3}$. A direct computation using 
(\ref{Akn})--(\ref{A0k}) and (\ref{An3a})--(\ref{An3c})
shows that these equations have only one common root
given by
\be\label{x32}
 \xi_{3} =\varkappa_q \,
 \frac{ \{2\}_q (5+3 \{2\}_q)}{ [2]_q [3]_q} \,,
\ee
which is the value corresponding to (\ref{serB}) for 
$s=\frac{3}{2}$. By Lemma~\ref{HR},
an R--matrix determined by (\ref{x12}) and (\ref{x32}) 
is unique and therefore it is the one given by~(\ref{serB}).
\qed\\[-0.5mm]

The proven proposition shows that, unlike the case of spin
$s=1$, only two out of four $sl_2$--invariant R--matrices 
(\ref{ser0a})--(\ref{ser0d}) have  
$U_q(sl_2)$--counterparts for spin $s=\frac{3}{2}$.
Actually, analysing the $n=3$ reduced YB equations, we
can extend this observation to higher spins as well.

\begin{prop}\label{Rs2}
Let $R(\lambda)$ be a $U_q(sl_2)$--invariant solution 
of the YB equation (\ref{YBE}) on $V_s^{\otimes 3}$
for a spin~$s\geq 2$ satisfying (\ref{gauge}) 
and~(\ref{ans2}). If 
$r_{2s-1}(\lambda) = 
\frac{[2s \+ \lambda]_q}{[2s \- \lambda]_q }$,
then 
\be\label{s22}
r_{2s-2}(\lambda) = 
\frac{[2s \+ \lambda]_q}{[2s \- \lambda]_q }
\frac{[2s\-1 \+ \lambda]_q}{[2s\-1 \- \lambda]_q } \,.
\ee
As a consequence, for $s\geq 2$, there exist no 
$U_q(sl_2)$--invariant regular R--matrices whose $q\to 1$ 
limit coincides with (\ref{ser0a}) or (\ref{ser0c}).
\end{prop}

\proof
Let $q=1+h$, $h \ll 1$, so that
$[t]_q = t + t(t\-1) h^2/3 + O(h^3)$. Since 
$A^{(s,n)}_{q}$ depends on $q$ smoothly, we  
have $A^{(s,n)}_{q} = A^{(s,n)}_{q=1} + O(h^2)$.
Consider the $n=3$ reduced YB equations (\ref{ddd})
where $\xi_0=0$, $\xi_1$ is as in (\ref{x1}),
and $\xi_2$ is given by the second expression in 
(\ref{sd3}). Using (\ref{Akn})--(\ref{A0k}) and 
(\ref{An3a})--(\ref{An3c}), we find the following 
$h$--expansions of these equations for 
$(a,b)=(0,1)$, $(0,2)$, and $(1,3)$, respectively
\be
\label{qha}
 0 &=& (5s^2 \- 3s) \,\xi_3^2 + (3 \- 6s) \,\xi_3 + 1 \\
\nonumber 
&& - \fr{2}{3} h^2 \, \bigl(  
 (25 s^4 \- 32 s^3 \+ 9 s^2)\,\xi_3^2 +
 (78 s^3 \- 81 s^2 \+ 21s) \,\xi_3
 \- 47 s^2 \+ 35 s \- 3 \bigr) + O(h^3) \,, \\ [0.5mm]
\label{qhb}
 0 &=& h^2 \, \bigl(
 (7s^2 \- 3s) \,\xi_3^2 + (3 \- 10s) \,\xi_3 + 3 
 \bigr) + O(h^3) \,, \\ [0.5mm]
\label{qhc}
 0 &=& (5s^2 \- 3s) \,\xi_3^2 + (3 \- 6s) \,\xi_3 + 1 \\
\nonumber 
&& + \fr{4}{3} h^2 \, \bigl(
 (19 s^4 \- 35 s^3 \+ 17 s^2 \- 3s)\,\xi_3^2 +
 (138 s^3 \- 201 s^2 \+ 96s \- 15) \,\xi_3
 \- 85 s^2 \+ 90 s \- 20
 \bigr) + O(h^3) \,.
\ee
We see that, in the zeroth order in $h$, (\ref{qhb}) is 
satisfied trivially, whilst (\ref{qha}) and (\ref{qhc}) yield 
the same quadratic equation which has the following roots
\be\label{x3}
 \xi_3 = \fr{1}{s} \,, \qquad 
 \xi_3 = \fr{1}{5s -3} \,.
\ee
Thus, for $q=1$, eqs.~(\ref{qha})--(\ref{qhc}) are
compatible (in particular, the first value in (\ref{x3}) 
corresponds to solutions of the type (\ref{ser0a})
and~(\ref{ser0c})).  

In the second order in $h$, (\ref{qhb}) has roots 
$\xi_3 = \frac{1}{s}$ and $\xi_3 = \frac{1}{7s -3}$.
But, for (\ref{qha}) and (\ref{qhc}), the $h^2$ corrections 
to the first value in (\ref{x3}) are 
\be
 \xi_3 = \fr{1}{s} + h^2 \, \bigl(
 \fr{28}{3}s \- 6 \bigr) + O(h^3) \,, \qquad
 \xi_3 = \fr{1}{s} + h^2 \, \bigl(
 \fr{46}{3} \- \fr{4}{s} \- 12s \bigr) + O(h^3) \,.
\ee
Therefore already in the second order in $h$ compatibility 
of (\ref{qha})--(\ref{qhc}) is lost. Which implies that
the second expression in (\ref{rn2}) for $r_{2s-2}(\lambda)$
is ruled out. And, as follows {} from the analysis of 
the previous subsection, the only possible form of 
$r_{2s-2}(\lambda)$ is the first expression in (\ref{rn2}).
An R--matrix with such spectral coefficient cannot be 
a~$q$--deformation of (\ref{ser0a}) or~(\ref{ser0c}) for
$s \geq 2$, because the corresponding value of $\xi_2$ 
does not vanish in the limit~$q \to 1$.
\qed\\[-0.5mm]

\rem 
Notice that the coefficient $r_{3}(\lambda)$ of the exceptional 
solution (\ref{s3a}) corresponds (after rescaling 
$\lambda\to\lambda/6$) to the second value in (\ref{x3}). 
As we have seen in the proof of Proposition~\ref{Rs2}, this value 
is not a root of (\ref{qhb}) for $h\neq 0$. Therefore, we conclude 
that (\ref{s3a}) has no regular $U_q(sl_2)$--invariant counterpart.

\section*{}
\vspace*{-7mm}
{\bf Acknowledgments:} 
This work was supported by the INTAS 
grant YS--03-55-962 and by the Russian Fund for Fundamental
Research grants 02-01-00085 and 03-01-00593.
I~am grateful to V.~Schomerus for his kind hospitality at the 
SPhT CEA--Saclay, where a part of this work was done.

\subsection*{Appendix A}

\begin{lem}\label{ru}
Let $R(\lambda)$ be a $U_q(sl_2)$--invariant solution of 
the YB equation (\ref{YBE}) on $V_s^{\otimes 3}$ satisfying 
conditions (\ref{reg}) and (\ref{gauge}). Then 
$R(\lambda)$ 
is unitary, i.e., it satisfies (\ref{uni}) as well.
\end{lem}

\proof Eq.~(\ref{ans1}) ensures that $R(\lambda)$ commutes 
with $R(\mu)$. Introduce $X^\lambda = R(\lambda)R(-\lambda)$.
Then the YB equation for $\mu=-\lambda$ implies that 
$X_{\one\two}^\lambda = X_{\two\three}^\lambda$.
Applying ${\rm\tr}_{\two\three}$ and 
${\rm\tr}_{\three}$ to this equality (along the 
lines of~\cite{Ke}, where a less trivial equation 
$X_{\one\two} - X_{\two\three}=Z_{\one\two\three}$
was considered), one infers that
$X^\lambda = c\, {\mathbb E}$, $c$ being a scalar constant.
On the other hand, we have $X^\lambda = P^{2s} + \ldots$,
according to~(\ref{gauge}). Hence $c=1$ and
$X^\lambda ={\mathbb E}$.
\qed

\subsection*{Appendix B}
The co-multiplication (\ref{del}) determines the structure
of the Clebsch--Gordan (CG) decomposition of tensor products
of irreducible representations. The corresponding CG
coefficients and 6--$j$ symbols were derived and 
studied in~\cite{CG}. The particular 6--$j$ symbol
which appeared in (\ref{A1}) is given by
\be
\label{6j}
 && \FRW{s}{s}{s}{3s \- n}{2s \- k}{2s \- k^\prime}_q =
 F^s_k \, F^s_{k'} \, \sum_{l} (-1)^l [l\+1]!  
 \, \Bigl(  [l\-4s\+k]! \, [l\-4s\+k^\prime]! \Bigm. \\
\nonumber
&& \qquad \times   \Bigm. 
 [l\-6s\+n\+k]! \,[l\-6s\+n\+k^\prime]! \,
 [6s\-n\-l]! \, [6s\-k\-k^\prime\-l]! \, 
 [8s\-n\-k\-k^\prime\-l]! \Bigr)^{-1} \,,
\ee
where
\be
\label{Fsk}
 F^s_{k} = [2s \- k]! \, \Bigl( \frac{[k]!\, 
 [n\-k]! \, [2s\-n\+k]! \, [4s\-n\-k]!}%
 {[4s\-k\+1]! \, [6s\-n\-k\+1]!} \Bigr)^{\frac{1}{2}}
\ee
and the $q$--factorial is defined as  
$[n]! =\prod_{k=1}^n [k]_q$ for $n\in {\mathbb Z}_+$ 
and~$[0]!=1$. The sum in (\ref{6j}) runs over
those $l$ for which the arguments of the $q$--factorials
are non--negative.

For $n\leq 2s$ and $k'\=0$ or $k'\=n$, the sum on the 
r.h.s.~of (\ref{6j}) contains only one term 
($l\=6s\-n$ or $l\=6s\-n\-k$, respectively), and we obtain 
\be 
\label{Akn}
  A^{(s,n)}_{k,n}  &=&
 \frac{(-1)^k \,\sqrt{[4s\-2k\+1]_q} }{[2s \- n]!} \\
\nonumber
 &\times& 
 \Bigl( \frac{[n]!\, [2s]! \, [2s\-n\+k]! \, [4s\-n\-k]!
 \, [4s\-2n\+1]! \, [6s\-n\-k\+1]! }%
 { [k]! \, [n\-k]! \, [4s\-k\+1]! \, [4s\-n\+1]! \,
 [6s\-2n\+1]!} \Bigr)^{\frac{1}{2}} ,\\ [1mm]
\label{A0k}
  A^{(s,n)}_{0,k}  &=&
 [2s]! \,\sqrt{[4s\-2k\+1]_q} \\ 
\nonumber
 &\times& \Bigl( \frac{[n]!\, [4s\-n]! \, [4s\-n\-k]!
 [6s\-n\+1]! }%
 { [k]! \, [n\-k]! \, [2s\-n\+k]! \, [2s\-n]! \,
 [4s]! \, [4s\-k\+1]! \, [6s\-n\-k\+1]! 
 } \Bigr)^{\frac{1}{2}} .
\ee

\subsection*{Appendix C}

{\em Proof} of Proposition~\ref{Aprop}:\par
$i)$ 
Applying the CG decomposition to the $\{23\}$ and $\{12\}$
components in (\ref{ek}) and (\ref{ek'}), respectively,
and using the orthonormality of the basis of $V_s$, 
$\langle p|p' \rangle = \delta_{p p'}$, it is straightforward 
to find that the scalar product in (\ref{defAn}) is given by
\be
\label{An}
 A^{(s,n)}_{k k^{\prime} } &=& \!\!\!
 \sum_{m, m^{\prime} }
 \FCG{s}{m}{2s \- k}{3s \- n \- m}{3s \- n}{3s \- n}_q \, 
 \FCG{2s \- k^{\prime}}{3s \- n \- m^{\prime}}%
  {s}{m}{3s \- n}{3s \- n}_q  \\
\nn
 & \times & \!\!\!
 \FCG{s}{3s \- n \- m \- m^{\prime}}{s}{m^{\prime}}%
    {2s \- k}{3s \- n \- m}_q \, 
 \FCG{s}{m}{s}{3s \- n \- m \- m^{\prime}}%
    {2s \- k^{\prime}}{3s \- n \- m^{\prime}}_q  \,.
\ee
In order to carry out the summation over $m$
we invoke the following identity~\cite{CG}
\be
\nonumber
&& \sum_{m} \, \FCG{a}{m}{b}{m' \- m}{e}{m'}_q \,
 \FCG{b}{m' \-m}{d}{m''}{f}{m'' \+ m' \- m}_q \,
 \FCG{a}{m}{f}{m'' \+ m' \- m}{c}{m'' \+ m' }_q \\
\label{CGRW}
&& \qquad = (-1)^{a\+b\+c\+d} \, \sqrt{[2e\+1]_q \,[2f\+1]_q} \, 
 \FCG{e}{m'}{d}{m''}{c}{m'' \+ m'}_q \,
 \FRW{a}{d}{b}{c}{e}{f}_q \,.
\ee
After this the summation over $m'$ reduces to
\be
\sum_{m'} \FCG{2s \- k'}{3s \- n \- m'}%
  {s}{m'}{3s \- n}{3s \- n}_q^2 = 
 \langle n; k' | n; k' \rangle = 1 \,.
\ee
The remaining factors in (\ref{CGRW}) yield the 
r.h.s.~of~(\ref{A1}).

$ii)$ The self--duality of the transition matrix 
with respect to $q \rightarrow q^{-1}$ follows {}from
the fact that 6--$j$ symbols are invariant with respect to
this operation (because, unlike the CG coefficients, 
they are expressed entirely in terms of 
$q$--numbers~\cite{CG}). 

$iii)$ The obvious invariance of (\ref{6j})
with respect to $k \leftrightarrow k'$ implies
that the transition matrix is symmetric. Since
$A^{(s,n)}_q$ is orthogonal by construction,
we conclude that $A^{(s,n)}_q$ coincides with 
its inverse. 

$iv)$ Formula (\ref{A6a}) is obvious from~(\ref{n1}).
The duality relation (\ref{A6b}) in terms of
matrix entries looks as follows
\be
\label{nn'}
 & A^{(s,n)}_{k,k'}  = 
  A^{(\tilde{s}, \tilde{n})}_{\tilde{k},\tilde{k}'}  \,, & \\
\label{ss'}
 & \tilde{s}=2s\-\frac{n}{2} \,, \quad
 \tilde{n}=6s\-2n \,, \quad
 \tilde{k}=k\-n\+2s \,, \quad \tilde{k}'=k'\-n\+2s \,, &
\ee
where $0 \leq k, k' \leq n$.
The shifts in $\tilde{k}$, $\tilde{k}'$ are necessary in
order to satisfy~(\ref{krange}) (notice that
$2s\-n=\tilde{n}\-2\tilde{s} \geq 0$). Eq.~(\ref{nn'})
is checked straightforwardly by making the change of 
variables (\ref{ss'}) in (\ref{A1}) and using the 
explicit expressions (\ref{6j})--(\ref{Fsk}). 
This completes the proof. \qed

\newcommand{\my}[6]{{#1:} {\em #2} {\bf #3} {(#4)} {#5}}
\small  \setlength{\itemsep}{-3pt}


\begin{thebibliography}{11}

\bibitem{KuR}
\my{P.P.~Kulish and N.Yu.~Reshetikhin}
{Zapiski Nauchn.\ Semin.\ LOMI}{101}{1981}{101}
\par\ (Engl.\ transl.\my{}{J.\ Sov.\ Math.}{23}{1983}{2435}\ );\\
\my{E.K.~Sklyanin}
{Funct.\ Analysis and Appl.}{16}{1982}{263}


\bibitem{SkDr}
\my{E.K.~Sklyanin}
{Uspekhi\ Mat.\ Nauk}{40}{1985}{214}\ ;\\
\my{M.~Jimbo}{Lett.\ Math.\ Phys.}{10}{1985}{63}\ ;\\
\my{V.G.~Drinfeld}
{Dokl.\ Akad.\ Nauk}{283}{1985}{1060}
\par\ (Engl.\ transl.\my{}{Sov.\ Math.\ Dokl.}{32}{1985}{254}\ )

\bibitem{Ro}
\my{M.~Rosso}
{Commun.\ Math.\ Phys.}{117}{1988}{581}

\bibitem{RM}
\my{J.B.~McGuire}
{J.\ Math.\ Phys.}{5}{1964}{622}\ ;\\
\my{C.N.~Yang}
{Phys.\ Rev.\ Lett.}{19}{1967}{1312}\ 

\bibitem{KRS}
\my{P.P.~Kulish, N.Yu.~Reshetikhin, and E.K.~Sklyanin}
{Lett. Math. Phys.}{5}{1981}{393}

\bibitem{ZZ}
\my{A.B.~Zamolodchikov and Al.B.~Zamolodchikov}
{Annals Phys.}{120}{1979}{253}

\bibitem{Ba}
\my{R.J.~Baxter}{J.~Stat.~Phys.}{28}{1982}{1}

\bibitem{Ke}
\my{T.~Kennedy}{J.~Phys.}{A25}{1992}{2809}

\bibitem{Jim} 
\my{M.~Jimbo}{Commun. Math. Phys.}{102}{1986}{537}

\bibitem{TL}
\my{H.N.V.~Temperley and E.H.~Lieb}
{Proc.\ Roy.\ Soc.\ Lond.}{A322}{1971}{251}

\bibitem{CG}
\my{A.N.~Kirillov and N.Yu.~Reshetikhin}
{Adv. Ser. Math. Phys.}{7}{1989}{285}\ ;\\
\my{M.~Nomura}{J.~Math.\ Phys.}{30}{1989}{2397}

\bibitem{KS}
\my{P.P.~Kulish and E.K.~Sklyanin}
{Lect.\ Notes Phys.}{151}{1982}{61}

\bibitem{Fad1} 
{L.D.~Faddeev:}  
{\em How algebraic Bethe ansatz works for integrable model}.
In: {\em Sym\'{e}tries quantiques, Les Houches 1995}
(North-Holland, 1998) 149 [hep-th/9605187]

\bibitem{Jo}
\my{V.~Jones}
{Commun.\ Math.\ Phys.}{125}{1989}{459}\ ;\\
Z.Q.~Ma: {\em Yang--Baxter equation and quantum enveloping 
algebras} (World Scientific, 1993)

\end{thebibliography}
\end{document}